\documentclass[twoside]{article}

\usepackage{amsmath,bbm,mor}
\usepackage{graphicx}
\usepackage{epstopdf}

\received{}
\revised{}
\pubyear{0x}
\pubmonth{Xxxxxxx}
\volume{xx}
\issue{x}
\pages{xxx--xxx}
\firstpage{0xxx}
\DOI{10.1287/moor.xxxx.xxxx}
\startpagenumber{1}

\title{Strategies for prediction under imperfect monitoring}
\ShortTitle{Strategies for prediction under imperfect monitoring}
\ShortAuthors{Lugosi, Mannor, and Stoltz}
\NumberOfAuthors{3}
\FirstAuthor{G\'abor Lugosi}
\FirstAuthorAddress{ICREA and Department of Economics, Pompeu Fabra University, Barcelona, Spain}
\FirstAuthorEmail{lugosi@upf.es}
\FirstAuthorURL{http://www.econ.upf.es/~lugosi}
\SecondAuthor{Shie Mannor}
\SecondAuthorAddress{Department of Electrical \& Computer Engineering, McGill University, Montreal, Canada}
\SecondAuthorEmail{shie.mannor@mcgill.ca}
\SecondAuthorURL{http://www.ece.mcgill.ca/~smanno1/}
\ThirdAuthor{Gilles Stoltz}
\ThirdAuthorAddress{
D\'epartement de Math\'ematiques et Applications,
Ecole Normale Sup\'erieure, CNRS, Paris, France \and
HEC Paris School of Management, CNRS, Jouy-en-Josas, France}
\ThirdAuthorEmail{gilles.stoltz@ens.fr}
\ThirdAuthorURL{http://www.dma.ens.fr/~stoltz}

\MSCcodes{
Primary:
91A20, 
62L12; 
Secondary:
68Q32, 
}
\ORMScodes{Primary:  computer science--artificial intelligence, decision analysis--sequential;
secondary:  games/group decisions--noncooperative}

\renewcommand\S{{\cal S}}
\def\cS{{\cal S}}
\def\cH{{\cal H}}

\def\cP{{\cal P}}

\def\F{{\cal F}}
\def\cF{{\cal F}}

\newcommand{\EXP}{\mathbbm{E}}
\newcommand{\R}{\mathbbm{R}}
\def\IND{{\mathbbm 1}}

\newcommand{\bp}{\boldsymbol{p}}
\newcommand{\bq}{\boldsymbol{q}}
\newcommand{\bu}{\boldsymbol{u}}
\newcommand{\bz}{\boldsymbol{z}}

\newcommand{\bw}{\boldsymbol{w}}

\newcommand{\be}{\boldsymbol{e}}
\newcommand{\by}{\boldsymbol{y}}
\newcommand{\bx}{\boldsymbol{x}}
\newcommand{\bb}{\boldsymbol{b}}

\renewcommand{\epsilon}{\varepsilon}

\newcommand{\wh}{\widehat}
\newcommand{\wt}{\widetilde}
\newcommand{\ol}{\overline}
\newcommand{\Var}{\mathrm{Var}}
\newcommand{\defeq}{\stackrel{\rm def}{=}}
\newcommand{\norm}[1][\cdot]{\ensuremath{\left\Arrowvert #1 \right\Arrowvert}}
\newcommand{\abs}[1][\cdot]{\ensuremath{\left| #1 \right|}}

\newlength{\minipagewidth}
\newcommand{\bookbox}[1]{
\par\medskip\noindent
\framebox[\textwidth]{
\begin{minipage}{\minipagewidth} {#1} \end{minipage} } \par\medskip }

\begin{document}
\maketitle

\bibliographystyle{plain}

\begin{abstract}
We propose simple randomized strategies for sequential
decision (or prediction)
under imperfect monitoring, that is, when the
decision maker (forecaster) does not
have access to the past outcomes but rather to a
feedback signal.
The proposed strategies are consistent in the sense that they
achieve, asymptotically, the best possible average reward
among all fixed actions.
It was Rustichini \cite{Rus99} who first proved the existence of such
consistent predictors. The forecasters presented here offer the
first constructive proof of consistency. Moreover, the proposed
algorithms are computationally efficient.
We also establish upper bounds for the rates
of convergence. In the case of deterministic
feedback signals,
these rates are optimal up to logarithmic terms.
\end{abstract}

\section{Introduction}

In sequential decision problems a decision maker (or forecaster)
tries to predict the outcome of a certain unknown process
at each (discrete) time instance and  takes an action accordingly.
Depending on the outcome of the predicted event and the action taken,
the decision maker receives a reward. Very often, probabilistic
modeling of the underlying process is difficult. For such situations
the prediction problem can be formalized as a repeated game between
the decision maker and the environment.
This formulation goes back to the 1950's when
Hannan~\cite{Han57}
and Blackwell~\cite{Bla56a} showed that the decision maker has a
randomized strategy that guarantees,  regardless of the outcome sequence,
an average asymptotic reward as
high as the maximal reward one could get by knowing the empirical
distribution of the outcome sequence in advance.
Such strategies are called {\em Hannan consistent}.
To prove this result,
Hannan and Blackwell assumed that the decision maker has full
access to the past outcomes.
This case is termed the {\em full information}
or the {\em perfect monitoring} case.
However, in many important applications,
the decision maker has limited information about the
past elements of the sequence to
be predicted. Various models of limited feedback
have been considered in the literature. Perhaps the best known
of them is the so-called {\em multi-armed bandit problem}
in which the forecaster is only informed
of its own reward but not the actual outcome;
see Ba\~nos~\cite{Ban68}, Megiddo~\cite{Meg80}, Foster and Vohra~\cite{FoVo98}, Auer, Cesa-Bianchi, Freund, and Schapire~\cite{ACFS02}, Hart and Mas Colell~\cite{HaMa99,HaMa02}.
For example, it is shown in \cite{ACFS02} that Hannan consistency
is achievable in this case as well.

Sequential decision problems like the ones considered in this paper have been studied in
different fields under various names such as repeated games, regret minimization,
on-line learning, prediction of individual sequences, and sequential prediction. The vocabulary of
different sub-communities differ. Ours is perhaps closest to that used by learning theorists.
For a general introduction and survey of the sequential prediction problem
we refer to Cesa-Bianchi and Lugosi \cite{CeLu06}.

In this paper we consider a general model in which the information
available to the forecaster is a general given  (possibly
randomized) function
of the outcome and the decision of the forecaster.
It is well understood under what conditions Hannan consistency
is achievable in this setup, see
Piccolboni and Schindelhauer \cite{PiSc01} and
Cesa-Bianchi, Lugosi, and Stoltz \cite{CeLuSt06}.
Roughly speaking, this is possible whenever, after suitable transformations
of the problem, the reward matrix can be expressed as a linear function of the
matrix of (expected) feedback signals. However, this condition is
not always satisfied and then the natural question is what
the best achievable performance for the decision maker is.
This question was answered by Rustichini \cite{Rus99} who characterized
the maximal achievable average reward that can be guaranteed asymptotically
for all possible outcome sequences (in an almost sure sense).

However, Rustichini's proof of achievability is not constructive.
It uses abstract {\em approachability} theorems due to
Mertens, Sorin, and Zamir \cite{MeSoZa94}
and it seems unlikely that
his proof method can give rise to computationally efficient
prediction algorithms, as noted in the conclusion of \cite{Rus99}.
A simplified efficient approachability-based strategy in the special case
where the feedback is a function of the action of nature alone was shown in Mannor and Shimkin \cite{MaSh03}.
In the general case, the simplified approachability-based strategy of \cite{MaSh03} falls short of the
maximal achievable average reward characterized by Rustuchini \cite{Rus99}.
The goal of this paper is to develop computationally efficient forecasters in the general prediction
problem under imperfect monitoring that
achieve the best possible asymptotic performance.

We introduce several forecasting strategies that exploit
some specific properties of the problem at hand. We separate
four cases, according to whether the feedback signal only depends on
the outcome or both on the outcome and the forecaster's action
and whether the feedback signal is deterministic or not. We design
different prediction algorithms for all four cases.

As a by-product, we also obtain finite-horizon performance
bounds with explicit guaranteed rates of convergence in terms of
the number $n$ of rounds the prediction game is played. In the case
of deterministic feedback signals these rates are optimal up to logarithmic
factors. In the random feedback signal case we do not know if it is
possible to construct forecasters with a significantly smaller regret.

A motivating example for such a prediction problem arises naturally in
multi-access channels that are prevalent in both wired and wireless
networks. In such networks, the communication medium is shared between
multiple decision makers. It is often technically difficult to
synchronize between the decision makers.  Channel sharing protocols,
and, in particular, several variants of spread spectrum, allow multiple
agents to use the same channel (or channels that may interfere with
each other) simultaneously.  More specifically, consider a wireless
system where multiple agents can choose in which channel to transmit
data at any given time. The quality of each channel may be different
and interference from other users using this channel (or other
``close'' channels) may affect the base-station reception.  The
transmitting agent may choose which channel to use and how much power
to spend on every transmission. The agent has a tradeoff between the
amount of power wasted on transmission and the cost of having its
message only partially received.  The transmitting agent may not
receive immediate feedback on how much data were received in the base
station (even if feedback is received, it often happens on a much
higher layer of the communication protocol). Instead, the transmitting
agent can monitor the transmissions of the other agents. However,
since the transmitting agent is physically far from the base-station
and the other agents, the information about the channels chosen by
other agents and the amount of power they used is imperfect. This
naturally abstracts to an online learning problem with imperfect
monitoring.

The paper is structured as follows.
In the next section we formalize the prediction problem we
investigate, introduce the  target quantity, that is,
the best achievable reward,
and the notion of regret.
In Section \ref{subgrad} we describe some analytical properties
of a key function
$\rho$, defined in Section \ref{setup}.
This function represents the worst possible average reward for a given vector of observations
and is needed in our analysis.
In Section \ref{aa} we consider the simplest special case
when the actions of the forecaster do not influence the
feedback signal, which is, moreover, deterministic. This case
is basically as easy as the full information case and
we obtain a regret bound of the order of $n^{-1/2}$
(with high probability) where $n$ is the number of rounds
of the prediction game.
In Section \ref{ba} we study random feedback signals but still
with the restriction that it is only determined by the outcome.
Here we are able to obtain a regret of the order of
$n^{-1/4}\sqrt{\log n}$.
The most general case is dealt with in Section \ref{bb}.
The forecaster introduced there has a regret of the order
of $n^{-1/5}\sqrt{\log n}$.
Finally, in Section \ref{ab} we show that this
may be improved to $O(n^{-1/3})$ in the case of deterministic
feedback signals, which is known to be optimal (see \cite{CeLuSt06}).

\section{Problem setup, notation}
\label{setup}

The randomized prediction problem is described as follows.
Consider a sequential decision problem in which a forecaster has to predict
an outcome that may be thought of as an action taken by the environment.

At each round, $t=1,2,\ldots,n$, the forecaster
chooses an action $i\in\{1,\ldots,N\}$ and the environment
chooses an action $j\in\{1,\ldots,M\}$ (which we also call an ``outcome'').
The forecaster's reward $r(i,j)$ is the value of
a reward function $r: \{1,\ldots,N\}\times\{1,\ldots,M\}\to[0,1]$.
Now suppose that, at the $t$-th round,
the forecaster chooses a
probability distribution
$\bp_t = (p_{1,t},\ldots,p_{N,t})$ over the set of actions,
and plays action $i$ with probability $p_{i,t}$.
We denote the forecaster's (random) action at time $t$ by $I_t$.
If the environment chooses action $J_t\in\{1,\ldots,M\}$,
then the reward of the forecaster is $r(I_t,J_t)$. The prediction problem is
defined as follows:

\bookbox{
\begin{center}
{\sc Randomized prediction with perfect monitoring}
\end{center}
\textbf{Parameters:} number $N$ of actions, cardinality $M$ of outcome space,
reward function $r$, number $n$ of game rounds.

\smallskip\noindent
For each round $t=1,2,\ldots,n,$

\smallskip\noindent
\begin{itemize}
\item[(1)]
the environment chooses the next outcome $J_t$;
\item[(2)]
the forecaster chooses $\bp_t$ and determines
the random action $I_t$, distributed according to $\bp_t$;
\item[(3)]
the environment reveals $J_t$;
\item[(4)]
the forecaster receives a reward $r(I_t,J_t)$.
\end{itemize}
}

Note in particular that the environment may react to the forecaster's strategy
by using a possibly randomized strategy. Below,
the probabilities of the considered events are taken with respect
to the forecaster's and the environment's randomized strategies.
The goal of the forecaster is to minimize the average regret
and to enforce that
\[
\limsup_{n \to \infty} \left(
 \max_{i=1,\ldots,N} \frac{1}{n}\sum_{t=1}^n r(i,J_t)
  - \frac{1}{n} \sum_{t=1}^n r(I_t,J_t) \right) \leq 0 \ \ \mbox{a.s.},
\]
that is, the per-round realized differences between the cumulative reward
of the best fixed strategy $i\in \{1,\ldots,N\}$, in hindsight,
and the reward of the forecaster, are asymptotically non positive.
Denoting by $r(\bp,j) = \sum_{i=1}^N  p_i r(i,j)$ the
linear extension of the reward function $r$, the Hoeffding-Azuma inequality
for sums of bounded martingale differences (see \cite{Hoe63}, \cite{Azu67}),
implies that for any $\delta \in (0,1)$,
with probability at least $1-\delta$,
\[
\frac{1}{n}\sum_{t=1}^n r(I_t,J_t) \ge \frac{1}{n}\sum_{t=1}^n r(\bp_t,J_t)
   - \sqrt{\frac{1}{2n}\ln \frac{1}{\delta}}~,
\]
so it suffices to study the average expected reward
$(1/n)\sum_{t=1}^n r(\bp_t,J_t)$. Hannan~\cite{Han57}
and Blackwell~\cite{Bla56a} were the first to show the existence of
a forecaster whose regret is $o(1)$ for all possible behaviors of the
opponent.
Here we mention  a simple yet powerful forecasting
strategy known as the {\em exponentially weighted average} forecaster.
This forecaster selects, at time $t$, an action $I_t$ according
to the probabilities
\[
   p_{i,t} = \frac{ \exp\left(\eta\sum_{s=1}^{t-1} r(i,J_s) \right)}
         {\sum_{k=1}^N \exp\left(\eta\sum_{s=1}^{t-1} r(k,J_s) \right)}~,
\qquad i=1,\ldots,N,
\]
where $\eta>0$ is a parameter of the forecaster.
One of the basic well-known results in the theory of prediction
of individual sequences states that the regret of the exponentially
weighted average forecaster is bounded as
\begin{equation}
\label{expertsthm}
\max_{i=1,\ldots,N} \frac{1}{n}\sum_{t=1}^n r(i,J_t)
  - \frac{1}{n} \sum_{t=1}^n r(\bp_t,J_t)
  \le \frac{\ln N}{n\eta} + \frac{\eta}{8}~.
\end{equation}
With the choice $\eta = \sqrt{8\ln N/n}$ the upper bound becomes
$\sqrt{\ln N/(2n)}$.
Different versions of this result have been proved
by Littlestone and Warmuth~\cite{LiWa94},
Vovk~\cite{Vov90,Vov98},
Cesa-Bianchi, Freund, Haussler,
Helmbold, Schapire, and Warmuth~\cite{CeFrHaHeScWa97},
Cesa-Bianchi~\cite{Ces99},
see also Cesa-Bianchi and Lugosi~\cite{CeLu99}.

In this paper we are concerned with problems in which the forecaster
does not have access neither to the outcomes $J_t$ nor
to the rewards $r(i,J_t)$. The information
available to the forecaster at each round is called the {\em feedback signal}.
These feedback signals may depend on the outcomes $J_t$ only or
on the action--outcome pairs $(I_t,J_t)$ and may
be deterministic or drawn at random.
In the simplest case when the feedback signal is deterministic,
the information available to the forecaster is $s_t = h(I_t,J_t)$, given by a
fixed (and known) deterministic feedback function
$h: \{1,\ldots,N\} \times \{1,\ldots,M\} \to \S$
where $\S$ is the finite set of signals.
In the most general case, the feedback signal is governed by a random feedback function of the form
$H: \{1,\ldots,N\} \times \{1,\ldots,M\} \to \cP(\S)$
where $\cP(\S)$ is the set of probability distributions over the signals.
The received feedback signal $s_t$ is then drawn at random according to
the probability distribution $H(I_t,J_t)$ by using an external independent randomization.

To make notation uniform throughout the paper, we identify a
deterministic feedback function $h: \{1,\ldots,N\} \times \{1,\ldots,M\} \to \S$
with the random feedback function
$H: \{1,\ldots,N\} \times \{1,\ldots,M\} \to \cP(\S)$
which, to each pair $(i,j)$, assigns $\delta_{h(i,j)}$
where $\delta_s$ is the probability distribution
concentrated on the single element $s \in \cS$.

The sequential prediction problem under imperfect monitoring
is formalized in Figure~\ref{FigGame}.
\begin{figure}
\bookbox{
\begin{center}
{\sc Randomized prediction under imperfect monitoring}
\end{center}
\textbf{Parameters:}
number $N$ of actions, number $M$ of outcomes,
reward function $r$, random feedback function $H$,
number $n$ of rounds.

\smallskip\noindent
For each round $t=1,2\ldots,n$,

\smallskip\noindent
\begin{enumerate}
\item
the environment chooses the next outcome $J_t\in\{1,\ldots,M\}$ without revealing it;
\item
the forecaster chooses a probability distribution $\bp_t$
over the set of $N$ actions and draws an action
$I_t\in\{1,\ldots,N\}$ according to this distribution;
\item
the forecaster receives reward $r(I_t,J_t)$ and each action $i$ gets
reward $r(i,J_t)$, but none of these values is revealed to the forecaster;
\item
a feedback signal $s_t$ drawn at random according to $H(I_t,J_t)$ is revealed to the forecaster.
\end{enumerate}
}
\caption{\label{FigGame}
The game of randomized prediction under imperfect monitoring}
\end{figure}

In many interesting situations
the feedback signal the forecaster receives is
independent of the forecaster's action and only depends on the outcome,
that is,
for all $j=1,\ldots,M$, $H(\cdot, j)$ is constant. In other words,
$H$ depends on the outcome $J_t$ but not on the forecaster's
action $I_t$.
We will see that the prediction problem becomes significantly simpler in
this special case.
To simplify notation in this case,
we write $H(J_t) = H(I_t,J_t)$ for the feedback signal at time $t$
($h(J_t) = h(I_t,J_t)$ in case of deterministic feedback signals).
This setting includes the full-information case (when the outcomes
$J_t$ are revealed) but also the case of noisy observations
(when a random variable with distribution depending only on $J_t$
is observed), see Weissman and Merhav~\cite{WeMe01},
Weissman, Merhav, and Somekh-Baruch~\cite{WeMeSo01}.

Next we describe a reasonable goal for the forecaster and define
the appropriate notion of consistency. To this end, we
introduce some notation.
If $\bp=(p_1,\ldots,p_N)$ and $\bq = (q_1,\ldots,q_M)$
are  probability distributions over $\{1,\ldots,N\}$
and $\{1,\ldots,M\}$, respectively, then, with a slight abuse of notation,
we write
\[
    r(\bp,\bq) = \sum_{i=1}^N \sum_{j=1}^M p_i q_j r(i,j)
\]
for the linear extension of the reward function $r$.
We also  extend linearly the random feedback function in its second argument:
for a probability distribution $\bq = (q_1,\ldots,q_M)$ over $\{1,\ldots,M\}$,
define the vector in $\cP(\S)$
\[
H(i,\bq) = \sum_{j=1}^M q_j H(i,j)~, \qquad i=1,\ldots,N.
\]
Denote by $\cal F$ the convex set of all $N$-vectors
$H(\cdot,\bq) = (H(1,\bq),\ldots,H(N,\bq))$ of probability distributions
obtained this way when $\bq$ varies.
(${\cal F} \subset \cP(\S)^N$
is the set of feasible distributions over the signals).
In the case when the feedback signals only depend on the outcome, all components of
this vector are equal and we denote their common value by $H(\bq)$.
We note that in the general case, the set $\cal F$ is the convex
hull of the $M$ vectors $H(\cdot,j)$. Therefore, performing a Euclidean projection on $\cal F$
can be done efficiently using quadratic programming.

To each probability distribution $\bp$ over $\{1,\ldots,N\}$
and probability distribution $\Delta \in \cal F$,
we may assign the quantity
\[
   \rho(\bp,\Delta) = \min_{\bq: H(\cdot,\bq) = \Delta} r(\bp,\bq)~,
\]
which is the reward guaranteed by the
mixed action $\bp$ of the forecaster against any distribution of the
outcomes that induces the given distribution of
feedback signals $\Delta$.
Note that $\rho\in [0,1]$ and that $\rho$ is concave in $\bp$ (since it is an infimum
of linear functions; since this infimum is taken on a convex set, the infimum is indeed a
minimum). Finally, $\rho$ is also convex in $\Delta$ as the condition defining the minimum
is linear in $\Delta$.

To define the goal of the forecaster,
let $\ol\bq_n$ denote the empirical distribution of the
outcomes $J_1,\ldots,J_n$ up to round $n$. This distribution may
be unknown to the forecaster since the forecaster observes
the signals rather than the outcomes.
The best the forecaster can hope for
is an average reward close to
$\max_{\bp} \rho(\bp,H(\cdot,\ol\bq_n))$.
Indeed, even if
$H(\cdot,\ol\bq_n)$ was known beforehand, the maximal expected reward for
the forecaster would be
$\max_{\bp} \rho(\bp,H(\cdot,\ol\bq_n))$, simply because without any additional
information the forecaster cannot hope to do better than
against the worst element which is equivalent to $\bq$ as far as the signals are concerned.

Based on this argument, the (per-round) regret $R_n$ is defined
as the average
difference between the obtained cumulative reward and the target
quantity described above, that is,
\[
R_n = \max_{\bp} \rho(\bp,H(\cdot,\ol\bq_n)) - \frac{1}{n} \sum_{t=1}^n r(I_t,J_t)~.
\]
Rustichini \cite{Rus99} proves the existence of a forecasting strategy
whose per-round regret is guaranteed to satisfy
$\limsup_{n\to \infty} R_n \le 0$ with probability one,
for all possible imperfect monitoring problems.

Rustichini's proof is not constructive but in several special
cases constructive and computationally efficient
prediction algorithms have been proposed.
Among the partial solutions proposed so far, we mention
Piccolboni and Schindelhauer \cite{PiSc01} and
Cesa-Bianchi, Lugosi, and Stoltz \cite{CeLuSt06} who study the
case when
\[
\max_{\bp} \rho(\bp,H(\cdot,\ol\bq_n))
= \max_{i=1,\ldots,N} r(i,\ol\bq_n)
= \max_{i=1,\ldots,N} \frac{1}{n} \sum_{t=1}^n r(i,J_t)~.
\]
In this case strategies with a vanishing per-round regret
are called {\em Hannan consistent}.
In such cases the feedback is sufficiently rich so that
one may achieve the same asymptotic reward as in the full information
case, although the rate of convergence may be slower.
This case turns out
to be considerably simpler to handle than the general problem and
computationally tractable explicit algorithms have been
derived. Also, it is shown in \cite{CeLuSt06} that in this case
it is possible to construct strategies whose regret decreases
at a rate of $n^{-1/3}$ (with high probability)
and that this rate of convergence cannot be
improved in general.
(Note that Hannan consistency is achievable, for example, in
the adversarial multi-armed bandit problem, see Remark~\ref{rkHC}
in the Appendix.)
Mannor and Shimkin \cite{MaSh03} construct an approachability-based algorithm with vanishing regret for the
special case where the feedback signals depend only on the outcome.
In addition, Mannor and Shimkin discuss the more general case of feedback signals that depend on
both the action and the outcome and provide an algorithm that attains a relaxed goal comparing to
the one attained in this work.

\begin{figure}[t]
\begin{center}
\includegraphics[width=9.0cm,height=5.95cm]{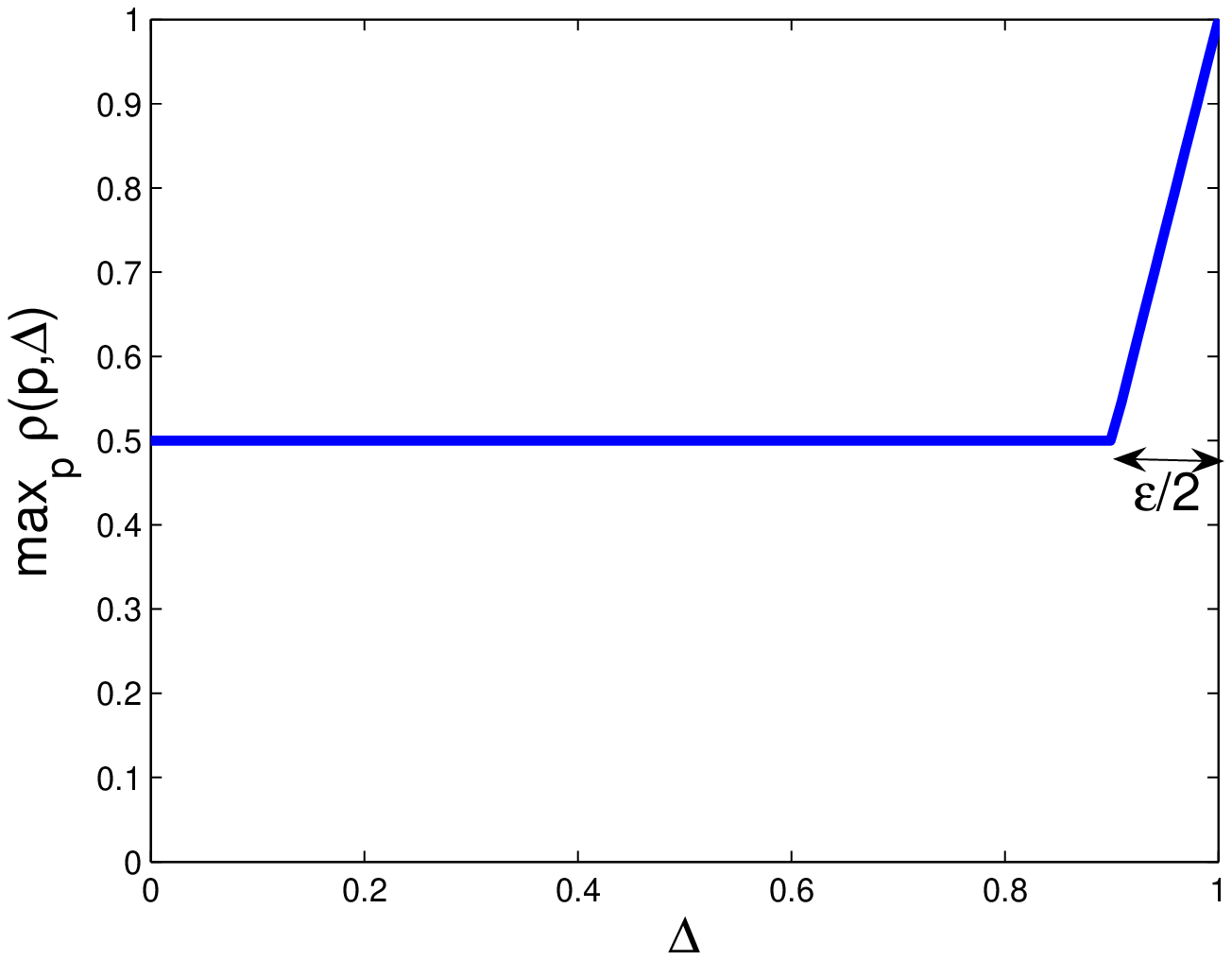}
\caption{The function $\Delta \mapsto
\max_{\bp} \rho(\bp,\Delta)$ for Example \ref{exmpl:2x3game}.
\label{fig:rhomax}}
\end{center}
\end{figure}

The following example demonstrates the
structure of the model.
\begin{example} \label{exmpl:2x3game}{\rm Consider the simple
game where $N=2$, $M=3$, $\cS=\{a,b\}$, and
the reward and feedback functions are as follows.
The reward function is described by the matrix
\[
\left[
\begin{array}{ccc}
  1 & 0 & 0 \\
  \frac{1}{2} & \frac{1}{2} & \frac{1}{2} \\
\end{array}
\right]
\]
To identify the possible distributions of the feedback signals we need to specify some elements of
$\cP(\cS)$. We describe such a member of $\cP(\cS)$ by the probability of observing $a$.
The feedback function is parameterized by some $\epsilon>0$ and is then given by
\[
\left[
\begin{array}{ccc}
  1 & 1-\epsilon & 0 \\
   1& 1-\epsilon & 0 \\
\end{array}
\right]~.
\]
In words, outcome $1$ leads to a deterministic feedback signal of $a$, outcome $3$
leads to a deterministic feedback signal of $b$, and outcome $2$
leads to a feedback signal of $a$ with probability $1-\epsilon$
and $b$ with probability $\epsilon$.
Note that the feedback signals depend only on the outcome and not on the action taken.
We recall that $\Delta$, as a member of $\cP(\cS)$,
is identified with the probability of observing the feedback signal $a$ and
it follows that $\cF$ is the interval $[0,1]$.
We now compute the function $\rho$.
Letting $p$ denote the probability of selecting the first
action (i.e., $\bp = (p,\,1-p)$), we have
\begin{eqnarray*}
\rho(\bp,\Delta) & = & \min_{\bq \,:\, q_1 + (1-\epsilon)q_2 = \Delta}
\left( p \, q_1 + (1-p) \frac{q_1 + q_2 + q_3}{2} \right)
=
\min_{\bq \,:\, \epsilon q_1 - (1-\epsilon) q_3 = \Delta - (1-\epsilon)} p \, q_1 + \frac{1-p}{2} \\
& = & \frac{1-p}{2} + \left\{
\begin{array}{ll}
\displaystyle{0}&
{\rm for}\ \Delta \le 1-\epsilon,
\vspace{5pt}\\ p \,\frac{\Delta - (1-\epsilon)}{\epsilon} & {\rm for}\ 1-\epsilon \le \Delta \le 1. \end{array} \right.
\end{eqnarray*}
Optimizing over $p$, we obtain
$$
\max_{\bp} \rho(\bp,\Delta) = \left\{
\begin{array}{ll}
\displaystyle{\frac{1}{2}}&
{\rm for}\ \Delta \le 1-\epsilon/2,
\vspace{5pt}\\ \frac{\Delta - (1-\epsilon)}{\epsilon} & {\rm for}\ 1-\epsilon/2 \le \Delta \le 1. \end{array} \right.
$$
The intuition here is that for $\Delta = 1$ there is certainty that the outcome is $1$ so that an action of $p=1$ is optimal.
For $\Delta \le 1-\epsilon$ the forecaster does not know if the outcome was consistently $2$ or some mixture of outcomes
$1$ and $3$.
By playing the second action, the forecaster can guarantee a reward of $1/2$. The function $\Delta \mapsto
\max_{\bp} \rho(\bp,\Delta)$ is
depicted in Figure \ref{fig:rhomax}.
}
\end{example}

In this paper we construct simple and computationally efficient
strategies whose regret vanishes with probability one.
The main idea behind the forecasters we introduce in the next sections
is based on the gradient-based strategies described, for example,
in Cesa-Bianchi and Lugosi \cite[Section 2.5]{CeLu06}.
Our forecasters use sub-gradients of concave functions.
In the next section we briefly recall some
basic facts on the existence, computation, and boundedness of these
sub-gradients.

\section{Some analytical properties of $\rho$}
\label{subgrad}

For a concave function $f$ defined over a convex subset of $\R^d$, a
vector $\bb(\bx) \in \R^d$ is a sub-gradient of $f$ at $\bx$ if
$f(\by)-f(\bx) \le \bb(\bx) \cdot (\by - \bx)$ for all $\by$ in the
domain of $f$. We denote  by $\partial f(\bx)$ the set of sub-gradients
of $f$ at $\bx$ which is also known as the sub-differential.  Sub-gradients always exist, that is, $\partial f(\bx)$
is non-empty in the interior of the domain of a concave function.  In
this paper, we are interested in sub-gradients of concave functions of
the form $f(\cdot) = \rho(\cdot,\wh{\Delta}_t)$, where $\wh\Delta_t$
is an observed or estimated distribution of feedback signal at round $t$. (For instance, in
Section~\ref{aa}, $\wh\Delta_t = \delta_{h(J_t)}$ is observed, in the
other sections, it will be estimated.)  In view of the exponentially
weighted update rules that are used below, we only evaluate these
functions in the interior of the definition domain (the simplex).
Thus, the existence of sub-gradients is ensured throughout.

In the general case, sub-gradients may be computed efficiently by
the simplex method.
However, their computation is often even simpler, as in the case
described in Section~\ref{aa}, that is, when one faces
deterministic feedback signals not depending on the actions of the forecaster.
Indeed, at round $t$, it is trivial whenever
$\bp \mapsto \rho(\bp,\delta_{h(J_t)})$ is
differentiable at the considered point $\bp_t$
since it is differentiable exactly
at those points at which it is locally linear, and thus
the gradient equals the column of the reward matrix
corresponding to the outcome $y_t$ for which
$r(\bp_t,y_t) = \rho(\bp_t,\delta_{h(J_t)})$.
But because $\rho(\cdot,\delta_{h(J_t)})$ is concave, the Lebesgue measure
of the set where it is non-differentiable equals zero.
It thus suffices to resort to the simplex method only at these points
to compute the sub-gradients.

Note that the components of the sub-gradients are always bounded by a constant
that depends on the game parameters.
This is the case since the $\rho(\cdot,\wh\Delta_t)$ are concave and
continuous on a compact set and are therefore Lipschitz, leading to a bounded sub-gradient.
In the sequel, we denote by $K$ the value
$\sup_{\bp} \sup_{\Delta}
\sup_{\bb \in \partial \rho(\bp,\Delta)}
\| \bb \|_{\infty}$
where $\partial \rho(\bp,\Delta)$ denotes the sub-gradient at $\bp$ of the
concave function $\rho(\cdot,\Delta)$ with $\Delta$ fixed.
This constant depends on the specific parameters of the game.
Since the parameters of the game are supposed to be known to the
forecaster, in principle, the forecaster can compute the value of $K$.
In any case, the value of $K$ can be bounded by
the supremum norm of the payoff function
as the following lemma asserts.

\begin{lemma}\label{le:Klessthan1}
The constant $K$ satisfies $K\leq 1$.
\end{lemma}

\begin{proof}
Fix $\Delta$ and consider $Z^\Delta = \{\bq:H(\cdot,\bq)=\Delta\}$.
Define $\varphi : (\bp,\bq) \in \R^n \times Z^\Delta \mapsto \varphi(\bp,\bq) \in \R$ as the linear
extension-restriction of $r$ to $\R^n\times Z^\Delta$, that is
$\varphi(\bp,\bq) = \sum_{i,j} p_i q_j r(i,j)$. Further, let
$Z^\Delta_0(\bp) = \{\overline{\bq}: \varphi(\bp,\overline{\bq})=\min_{\bq\in Z^\Delta}\varphi(\bp,\bq)\}$.
It follows that under our notation, for any probability distribution $\bp$, one has
$\rho(\bp,\Delta) = \min_{\bq\in Z^\Delta} \varphi(\bp,\bq)$. Now, from Danskin's theorem
(see, e.g., Bertsekas \cite{Ber95})
we have that the sub-differential satisfies
\[
\partial \rho(\bp,\Delta) = {\rm conv} \left( \frac{\partial \varphi(\bp,\bz)}{\partial \bp} : \bz\in Z^\Delta_0(\bp)\right)
\]
where ${\rm conv}(A)$ denotes the convex hull of a set $A$.
Since $r(i,j)\in [0,1]$, it follows that
$\| \partial \rho(\bp,\bz)/\partial \bp\|_\infty \leq 1$ for all $\bz\in Z^\Delta$.
Since the convex hull does not increase the infinity norm, the result follows.
\end{proof}

\begin{remark} {\rm The constant $K$ for the game described in Example \ref{exmpl:2x3game} is $1/2$.
However, the gradient of the function $\max_{\bp} \rho(\bp,\Delta)$ as a function
of $\Delta$ is $1/\epsilon$.
This happens because the $\bp$ that attains the maximum changes rapidly in the interval $[1-\epsilon/2,1]$.
We further note that $K$ may be much smaller than 1. Since our regret bounds below depend on
$K$ linearly, having a tighter bound on $K$ can lead to considerable convergence rate speedup;
see Remark \ref{r:eta}.}
\end{remark}

\section{Deterministic feedback signals only depending on outcome}
\label{aa}

We start with the simplest case when the feedback signal is
deterministic and does not depend on the action $I_t$ of the forecaster.
In other words, after making the prediction at time $t$,
the forecaster observes $h(J_t)$.
This simplifying assumption may be naturally satisfied in
applications in which the forecaster's decisions do not effect
the environment.

In this case, we group the outcomes according
to the deterministic feedback signal they are associated to. Each signal
$s$ is uniquely associated to a group of outcomes.
This situation is very similar to the case of full monitoring
except that rewards are measured
by $\rho$ and not by $r$.
This does not pose a problem since $r$ is lower bounded
by $\rho$ in the sense that for all $\bp$ and $j$,
\[
r(\bp,j) \geq \rho(\bp,\delta_{h(j)})~.
\]
As mentioned in the previous section, we introduce a
forecaster based on the sub-gradients of $\rho(\cdot,\delta_{h(J_t)})$,  $t = 1,2,\ldots$.
The forecaster requires a tuning parameter $\eta > 0$. The $i$-th component of $\bp_t$
is
\[
   p_{i,t}
= \frac{e^{\eta\sum_{s=1}^{t-1} \left( \wt{r} (\bp_s,\delta_{h(J_s)}) \right)_i }}
   {\sum_{j=1}^N e^{\eta\sum_{s=1}^{t-1} \left( \wt{r} (\bp_s,\delta_{h(J_s)}) \right)_j }}~,
\]
where $\left( \wt{r} (\bp_s,\delta_{h(J_s)}) \right)_i$ is the $i$-th component
of any sub-gradient $\wt{r}(\bp_s,\delta_{h(J_s)}) \in \partial \rho(\bp_s,\delta_{h(J_s)})$
of the concave function $\rho(\cdot,\delta_{h(J_s)})$.
This forecaster is inspired by a gradient-based predictor introduced
by Kivinen and Warmuth \cite{KW97}.

The regret is bounded as follows.
Note that the following bound and
the considered forecaster coincide with those of (\ref{expertsthm})
in case of perfect monitoring. (In that case, $\rho(\cdot, \delta_{h(j)}) = r(\cdot,j)$, the sub-gradients are given by $r$.)

\begin{proposition}
\label{prop:aa}
For all $\eta > 0$, for all strategies of the environment,
for all $\delta > 0$,
the above strategy of the forecaster ensures that,
with probability at least $1 - \delta$,
\[
R_n
\leq
\frac{\ln N}{\eta n} + \frac{K^2\eta}{2}
   + \sqrt{\frac{1}{2n} \ln \frac{1}{\delta}}~,
\]
where
$K$ is the bound on the sub-gradients considered above.
In particular, choosing $\eta \sim \sqrt{ (\ln N) / n}$ yields
$R_n = O(n^{-1/2} \sqrt{\ln (N/\delta)})$.
\end{proposition}

\begin{remark}
\label{r:eta}
The optimal choice of $\eta$ in the upper bound is
$K \sqrt{2 (\ln N)/n}$, which depends on the parameters $K$
and $n$. While the bound $K\le 1$ is available, this bound might be loose.
Sometimes the forecaster does not necessarily know
in advance the number
of prediction rounds and/or the value of $K$ may be difficult to compute.
In such cases one may
estimate on-line both the number of time rounds and $K$, using the
techniques of Auer, Cesa-Bianchi, and Gentile~\cite{AuCeGe02} and
Cesa-Bianchi, Mansour, and Stoltz~\cite{CeMaSt07} as follows. Writing
\[
K_t = \max_{s \leq t-1} \| \wt{r}(\bp_s,\delta_{h(J_s)}) \|_\infty~,
\]
and introducing a round-dependent choice of the tuning parameter $\eta
= \eta_t = C K_t \sqrt{(\ln N)/t}$ for a properly chosen constant $C$,
one may prove a regret bound that is a constant multiple of $K_n
\sqrt{(\ln N/\delta)/n}$ (that hold with probability at least
$1-\delta$).
Since the proof of this is a straightforward combination of the techniques
of the above-mentioned papers and our proof, the details are omitted.
\end{remark}

\begin{proof}
Note that since the feedback signals are deterministic,
$H(\ol\bq_n)$ takes the simple form
$H(\ol\bq_n) = \frac{1}{n} \sum_{t=1}^n \delta_{h(J_t)}$.
Now, for any $\bp$,
\begin{eqnarray*}
\lefteqn{
n \rho(\bp,H(\ol\bq_n)) - \sum_{t=1}^n r(\bp_t,J_t)   } \\*
& \le &
 n \rho(\bp,H(\ol\bq_n)) - \sum_{t=1}^n \rho(\bp_t,\delta_{h(J_t)})
\quad \mbox{(by the lower bound on $r$ in terms of $\rho$)} \\
& \le &
\sum_{t=1}^n \left( \rho(\bp,\delta_{h(J_t)}) - \rho(\bp_t,\delta_{h(J_t)}) \right)
\quad \mbox{(by convexity of $\rho$ in the second argument)}  \\
& \le &
\sum_{t=1}^n \wt{r}(\bp_t,\delta_{h(J_t)}) \cdot (\bp - \bp_t)
\quad \mbox{(by concavity of $\rho$ in the first argument)} \\
& \le &
\frac{\ln N}{\eta} + \frac{nK^2\eta}{2}
\quad
\mbox{(by (\ref{expertsthm}), after proper rescaling),}
\end{eqnarray*}
where at the last step we used the fact that the forecaster is just
the exponentially weighted average predictor based on the rewards
$(\wt{r}(\bp_s,\delta_{h(J_s)}))_i$ and that all these reward vectors have
components between $-K$ and $K$. The proof is concluded by the
Hoeffding-Azuma inequality, which ensures that, with probability at
least $1- \delta$,
\begin{equation}
\label{HAzconcludes}
\sum_{t=1}^n r(I_t,J_t) \geq \sum_{t=1}^n r(\bp_t,J_t) - \sqrt{\frac{n}{2} \ln \frac{1}{\delta}}~.
\end{equation}
\end{proof}

\section{Random feedback signals depending only on the outcome}
\label{ba}

Next we consider the case when the feedback signals do not depend
on the forecaster's actions, but, at time $t$, the signal $s_t$
is drawn at random
according to the distribution $H(J_t)$.
In this case the forecaster does not have a direct access to
\[
H(\ol\bq_n) = \frac{1}{n} \sum_{t=1}^n H(J_t)
\]
anymore, but only observes the realizations $s_t$ drawn at random
according to $H(J_t)$.
In order to overcome this problem, we
group together several consecutive time rounds
(say, $m$ of them) and estimate
the probability distributions according to which the signals have been drawn.

To this end, denote by $\Pi$ the Euclidean projection
onto $\F$
(since the feedback signals depend only on the outcome we may now view
the set $\F$ of feasible distributions over the
signals as a subset of $\cP(\cS)$, the latter being identified with a
subset of $\R^{|\cS|}$ in a natural way).
Let $m$, $1\le m\le n$, be a parameter of the algorithm.
For $b = 0, 1, \ldots$, we denote
\begin{equation}
\label{defPD-group}
   \wh\Delta^b = \Pi \left( \frac{1}{m} \sum_{t=bm+1}^{(b+1)m} \delta_{s_t} \right)~.
\end{equation}
For the sake of the analysis, we also introduce
\[
   \Delta^b = \frac{1}{m} \sum_{t=bm+1}^{(b+1)m} H(J_t)~.
\]
The proposed strategy is described in Figure \ref{figba}.  Observe
that the practical implementation of the forecaster only requires the
computation of (sub)gradients and of $\ell_2$ projections, which can
be done in polynomial time.
\begin{figure}[t]
\bookbox{
\textbf{Parameters:} Integer $m \geq 1$, real number $\eta > 0$.
\\
\textbf{Initialization:} $\bw^0 = (1,\ldots,1)$.

\smallskip\noindent
For each round $t=1,2,\ldots$

\smallskip\noindent
\begin{enumerate}
\item If  $bm+1 \leq t < (b+1)m$ for some integer $b$,
choose the distribution $\bp_t = \bp^b$ given by
\[
p_{k,t} = p_k^b = \frac{w^b_k}{\sum_{j=1}^N w^b_j}
\]
and draw an action $I_t$ from $\{1,\ldots,N\}$ according to it;
\item
if $t = (b+1)m$ for some integer $b$, perform the update
\[
    w^{b+1}_{k} = w^b_{k}\,e^{\eta\,\left(\wt{r}\left(\bp^b,\wh{\Delta}^b\right)\right)_k}
\qquad
    \mbox{\rm for each $k=1,\ldots,N$,}
\]
where for all $\Delta$, \
$\wt{r}(\cdot,\Delta)$ is a sub-gradient of $\rho(\cdot, \Delta)$
and $\wh\Delta^b$ is defined in (\ref{defPD-group}).
\end{enumerate}
}
\caption{\label{figba}
The forecaster for random feedback signals depending only on the outcome.
}
\end{figure}
The next theorem bounds the regret of the strategy which is of the
order of $n^{-1/4}\sqrt{\log n}$. The price we pay for
having to estimate the distribution is thus a deteriorated
rate of convergence (from the $O(n^{-1/2})$ obtained in the case
of deterministic feedback signals). We do not know whether this rate can be
improved significantly as we do not know of any nontrivial lower bound
in this case.

\begin{theorem}
\label{Th:RandomOutcomeOnly}
For all integers $m \geq 1$, for all $\eta > 0$, and for all $\delta > 0$,
the regret for any strategy of the environment is bounded,
with probability at least $1 - (n/m + 1) \delta$, by
\[
R_n \le
2\sqrt{2} \, L \, \frac{1}{\sqrt{m}} \sqrt{\ln \frac{2}{\delta}} + \frac{m \ln N}{n \eta}
+ \frac{K^2\eta}{2} + \frac{m}{n} + \sqrt{\frac{1}{2n} \ln \frac{1}{\delta}} ~,
\]
where $K\le 1$ and $L$ are constants that depend only on the parameters of the game.
The choices $m = \lceil \sqrt{n} \rceil$ and
$\eta \sim \sqrt{(m \ln N)/n}$ imply
$R_n = O(n^{-1/4} \sqrt{\ln (nN/\delta)})$
with probability of at least $1-\delta$.
\end{theorem}

\begin{remark}
Here again, $K$ and $L$ may, in principle, be computed or bounded (see
Lemma \ref{le:Klessthan1} and Remark~\ref{rem:boundL}) by the
forecaster.
If the horizon $n$ is known in advance (as it is assumed in this paper),
the values of $\eta$ and $m$ may be chosen to
optimize the upper bound for the regret.  Observe that while one
always have $K\le 1$, the value of $L$ (i.e., the Lipschitz constant
of $\rho$ in its second argument) can be arbitrarily large, see
Example \ref{exmpl:2x3game}.
If the horizon $n$ is unknown at the start of the game, the situation
is not as simple as in Section \ref{aa} (see Remark~\ref{r:eta}),
because now a time-dependent choice of
$\eta$ needs to be accompanied by an adaptive choice of the parameter
$m$ as well. A simple, though not very attractive, solution is the so-called
``doubling trick'' (see, e.g., \cite[p.17]{CeLu06}). According to this
solution, time is divided into periods of exponentially growing length
and in each period the forecaster is used as if the horizon was the
length of the actual period. At the end of each period the forecaster
is reset and started again with new parameter values. It is easy to see
that this forecaster achieves the same regret bounds, up to a constant
multiplier. We believe that a smoother solution should also work
(as in Remark \ref{r:eta}).
Since this seems like a technical endeavor
we do not pursue this issue further.
\end{remark}

\begin{proof}
We start by grouping time rounds $m$ by $m$.
For simplicity, we assume that $n = (B+1) m$ for some integer $B$;
if this is not the case, we consider the lower integer part of $n$
and bound the regret suffered in the last at most $m-1$ rounds
by $m$ (this accounts for the $m/n$ term in the bound).
For all $\bp$,
\begin{eqnarray*}
n \, \rho(\bp,H(\ol\bq_n)) - \sum_{t=1}^n r(\bp_t,J_t)
& = &
n \, \rho(\bp,H(\ol\bq_n))
- \sum_{b = 0}^B m \, r\left( \bp^b, \frac{1}{m} \sum_{t=bm+1}^{(b+1)m} \delta_{J_t} \right)
\\
& \leq &
\sum_{b = 0}^B \left( m \, \rho \left( \bp, \Delta^b \right)
- m \, r\left( \bp^b, \frac{1}{m} \sum_{t=bm+1}^{(b+1)m} \delta_{J_t} \right)
 \right) \\
& \leq &
m \sum_{b = 0}^B
\left( \rho \left( \bp, \Delta^b \right)
- \rho \left( \bp^b, \Delta^b \right)
 \right)~,
\end{eqnarray*}
where we used the definition of the algorithm,
convexity of $\rho$ in its second argument,
and finally, the definition of $\rho$ as a minimum.
We proceed by estimating  $\Delta^b$ by $\wh\Delta^b$.
By a version
of the Hoeffding-Azuma inequality for sums of
Hilbert space-valued
martingale differences proved by
Chen and White \cite[Lemma 3.2]{ChWh96}, and since the $\ell_2$ projection
can only help, for all $b$, with probability at least $1 - \delta$,
\[
\norm[\Delta^b - \wh\Delta^b]_2 \leq \sqrt{\frac{2 \ln \frac{2}{\delta}}{m}}~.
\]
By Proposition \ref{prop:uniformLip}, $\rho$  is uniformly Lipschitz
in its second argument (with constant $L$), and therefore
we may further bound as
follows. With probability $1-(B+1)\delta$,
\begin{eqnarray*}
m \sum_{b = 0}^B \left(
\rho \left( \bp, \Delta^b \right)
- \rho \left( \bp^b, \Delta^b \right)
 \right)
& \leq &
m  \sum_{b = 0}^B \left(
\rho \left( \bp, \wh\Delta^b \right) -
\rho \left( \bp^b, \wh\Delta^b \right)
+ 2 \, L \sqrt{\frac{2 \ln \frac{2}{\delta}}{m}} \right) \\
& = &
m \sum_{b = 0}^B \left(
\rho \left( \bp, \wh\Delta^b \right)
- \rho \left( \bp^b, \wh\Delta^b \right) \right)
+ 2 \, L (B+1) \sqrt{2 m \ln \frac{2}{\delta}}~.
\end{eqnarray*}
The term containing $(B+1) \sqrt{m} = n/\sqrt{m}$ is the first term in the upper bound.
The remaining part is bounded by using the same slope inequality argument
as in the previous section
(recall that $\wt{r}$ denotes a sub-gradient),
\begin{eqnarray*}
m \sum_{b = 0}^B \left(
\rho \left( \bp, \wh\Delta^b \right)
- \rho \left( \bp^b, \wh\Delta^b \right)
\right)
& \leq &
m \sum_{b = 0}^B \wt{r} \left( \bp^b, \wh\Delta^b \right) \, \cdot \, \left( \bp - \bp^b \right)  \\*
& \le &
m \left( \frac{\ln N}{\eta} + \frac{(B+1)K^2 \eta}{2} \right) = \frac{m \ln N}{\eta} + \frac{n K^2\eta}{2}
\end{eqnarray*}
where we used Theorem \ref{expertsthm} and the boundedness of the
function $\wt{r}$ between $-K$ and $K$.
The proof is concluded by the Hoeffding-Azuma inequality which,
as in (\ref{HAzconcludes}), gives the final term in the bound.
The union bound indicates that the obtained bound holds with
probability at least $1 - (B+2) \delta \geq 1 - (n/m + 1) \delta$.
\end{proof}

\section{Random feedback signals depending on action--outcome pair}
\label{bb}

We now turn to the general case, where the feedback signals are random and
depend on the action--outcome pairs $(I_t,J_t)$. The key is, again,
to exhibit efficient estimators of the (unobserved) $H(\cdot,\ol\bq_n)$.

Denote by $\Pi$ the projection, in the
Euclidian distance, onto $\F$
(where $\F$, as a subset of $(\cP(\cS))^N$, is identified with
a subset of $\R^{|\cS| \, N}$). For $b = 0, 1, \ldots$, denote
\begin{equation}
\label{defPD-group2}
   \wh\Delta^b = \Pi \left( \frac{1}{m} \sum_{t=bm+1}^{(b+1)m} \left[ \wh{h}_{i,t} \right]_{i=1,\ldots,N} \right)
\end{equation}
where the distribution $H(i,J_t)$ of the random signal $s_t$
received by action $i$ at round $t$ is estimated by
\[
\wh{h}_{i,t} = \frac{\delta_{s_t}}{p_{i,t}} \IND_{I_t = i}~.
\]
(This form of estimators is reminiscent of those
presented, e.g., in \cite{ACFS02, PiSc01, CeLuSt06}.)
We prove that the $\wh{h}_{i,t}$ are conditionally unbiased
estimators. Denote by $\EXP_t$ the conditional expectation with
respect to the information available to the forecaster at the
beginning of round $t$. This conditioning fixes the values of $\bp_t$
and $J_t$.
Thus,
\[
\EXP_t \left[ \wh{h}_{i,t} \right] = \frac{1}{p_{i,t}} \, \EXP_t \left[ \delta_{s_t} \IND_{I_t = i} \right] =
\frac{1}{p_{i,t}} \, \EXP_t \left[ H(I_t,J_t) \IND_{I_t = i} \right] = \frac{1}{p_{i,t}} H(i,J_t) p_{i,t} = H(i,J_t)~.
\]
For the sake of the analysis, introduce
\[
   \Delta^b = \frac{1}{m} \sum_{t=bm+1}^{(b+1)m} H(\cdot,J_t)~.
\]
The proposed forecasting strategy is described in Figure \ref{figbb}.
The mixing with the uniform distribution is needed, similarly to the forecasters
presented in \cite{ACFS02, PiSc01, CeLuSt06}, to ensure sufficient exploration
of all actions. Mathematically, such a mixing lower bounds the probability
of pulling each action, which will turn to be crucial in the proof of
Theorem~\ref{Th:RandomActionOutcome}.
\begin{figure}[t]
\bookbox{
\textbf{Parameters:} Integer $m \geq 1$, real numbers $\eta, \gamma > 0$.
\\
\textbf{Initialization:} $\bw^0 = (1,\ldots,1)$.

\smallskip\noindent
For each round $t=1,2,\ldots$

\smallskip\noindent
\begin{enumerate}
\item if  $bm+1 \leq t < (b+1)m$ for some integer $b$,
choose the distribution $\bp_t = \bp^b = (1-\gamma) \wt\bp^b + \gamma \bu$,
where $\wt\bp^b$ is defined component-wise as
\[
\wt{p}_k^b = \frac{w^b_k}{\sum_{j=1}^N w^b_j}
\]
and $\bu$ denotes the uniform distribution, $\bu = (1/N,\ldots,1/N)$;
\item draw an action $I_t$ from $\{1,\ldots,N\}$ according to it;
\item
if $t = (b+1)m$ for some integer $b$, perform the update
\[
    w^{b+1}_{k} = w^b_{k}\,e^{\eta\,\left(\wt{r}\left(\bp^b,\wh{\Delta}^b\right)\right)_k}
\qquad
    \mbox{\rm for each $k=1,\ldots,N$,}
\]
where for all $\Delta \in \cF$,
$\wt{r}(\cdot,\Delta)$ is a sub-gradient of $\rho(\cdot,\Delta)$
and $\wh\Delta^b$ is defined in (\ref{defPD-group2}).
\end{enumerate}
}
\caption{\label{figbb}
The forecaster for random feedback signals depending on action--outcome pair.
}
\end{figure}

Here again, the practical implementation of the forecaster
only requires the computation of (sub)gradients and of $\ell_2$ projections,
which can be done efficiently. The next theorem
states that the regret in this most general case is
at most of the order of $n^{-1/5}\sqrt{\log n}$.
Again, we do not know whether this bound can be improved significantly.
We recall that $K$ denotes an upper bound on the infinity norm of the sub-gradients (see
Lemma~\ref{le:Klessthan1}).
The issues concerning the tuning of the parameters considered in the
following theorem are similar to those discussed after the statement
of Theorem~\ref{Th:RandomOutcomeOnly}; in particular, the simplest
way of being adaptive in all parameters is to use the ``doubling trick".

\begin{theorem}
\label{Th:RandomActionOutcome}
For all integers $m \geq 1$, for all $\eta > 0$, $\gamma \in (0,1)$, and $\delta > 0$,
the regret for any strategy of the environment
is bounded, with probability at least $1 - (n/m+1) \delta$, as
\begin{eqnarray*}
R_n & \le &
2 L \, N \sqrt{\frac{2\abs[\cS]}{\gamma m} \ln \frac{2N \abs[\cS]}{\delta}}
+ 2 L \frac{N^{3/2} \sqrt{\abs[\cS]}}{3 \gamma m} \ln \frac{2N \abs[\cS]}{\delta} \\*
& & \qquad  + \frac{m \ln N}{n \eta} + \frac{K^2\eta}{2} + 2 K\,\gamma
+ \frac{m}{n} + \sqrt{\frac{1}{2n} \ln \frac{1}{\delta}}~,
\end{eqnarray*}
where $L$ and $K \leq 1$  are constants that depend on the parameters of the game.
The choices $m = \lceil n^{3/5} \rceil$, $\eta \sim \sqrt{(m \ln N)/n}$,
and $\gamma \sim n^{-1/5}$
ensure that,  with probability at least $1 - \delta$,
$R_n = O\left(n^{-1/5} N \sqrt{\ln \frac{Nn}{\delta}}
+ n^{-2/5} N^{3/2} \ln \frac{Nn}{\delta} \right)$.
\end{theorem}

\begin{proof}
The proof is similar to the one of Theorem~\ref{Th:RandomOutcomeOnly}.
A difference is that we bound the accuracy of the estimation of
the $\Delta^b$ via a martingale analog of Bernstein's inequality
due to Freedman \cite{Fre75}
rather than the Hoeffding-Azuma inequality.
Also, the mixing with the uniform distribution
in the first step
of the definition of the forecaster in Figure~\ref{figbb} needs to be handled.

We start by grouping time rounds $m$ by $m$.  Assume,
for simplicity, that $n = (B+1) m$ for some integer $B$ (this accounts, again, for the
$m/n$ term in the bound). As before, we get that,
for all $\bp$,
\begin{equation}\label{eq:nrhoboundedbymb}
 n \, \rho(\bp,H(\cdot,\ol\bq_n))
- \sum_{t=1}^n r(\bp_t,J_t)
\leq
m \sum_{b = 0}^B
\left( \rho \left( \bp, \Delta^b \right)
- \rho \left( \bp^b, \Delta^b \right)\right)
\end{equation}
and proceed by estimating $\Delta^b$ by $\wh\Delta^b$.
Freedman's inequality \cite{Fre75}
(see, also, \cite[Lemma A.1]{CeLuSt06}) implies that
for all $b = 0, 1, \ldots, B$, $i = 1,\ldots,N$, $s \in \cS$, and $\delta > 0$,
\[
\abs[\Delta^b_i(s) - \frac{1}{m} \sum_{t=bm+1}^{(b+1)m} \wh{h}_{i,t}(s)]
\leq \sqrt{2 \frac{N}{\gamma m} \ln \frac{2}{\delta}} + \frac{1}{3} \frac{N}{\gamma m} \ln \frac{2}{\delta}
\]
where $\wh{h}_{i,t}(s)$ is the probability mass put on $s$ by
$\wh{h}_{i,t}$ and $\Delta^b_i(s)$ is the $i$-th component of $\Delta^b$.
This is because the sums of the conditional variances are bounded as
\[
\sum_{t=bm+1}^{(b+1)m}
  \Var_t \left(\frac{\IND_{I_t = i, s_t = s}}{p_{i,t}}\right)
 \leq
\sum_{t=bm+1}^{(b+1)m} \frac{1}{p_{i,t}} \leq \frac{m N}{\gamma}
\]
where the second inequality follows from the lower bound $\gamma/N$
on the components of $\bp_t$ (ensured by the mixing step in the definition
of the forecaster).
Summing (since the $\ell_2$ projection can only help), the union bound
shows that for all $b$, with probability at least $1 - \delta$,
\[
\norm[\Delta^b - \wh\Delta^b]_2 \leq d \defeq \sqrt{N \abs[\cS]} \left( \sqrt{2 \frac{N}{\gamma m} \ln \frac{2N \abs[\cS]}{\delta}}
+ \frac{1}{3} \frac{N}{\gamma m} \ln \frac{2N \abs[\cS]}{\delta} \right)~.
\]
By using uniform Lipschitzness of $\rho$ in its second argument
(with constant $L$; see Proposition~\ref{prop:uniformLip}), we may further
bound (\ref{eq:nrhoboundedbymb}) with probability $1-(B+1)\delta$ by
\begin{eqnarray*}
m \sum_{b = 0}^B \left(
\rho \left( \bp, \Delta^b \right)
- \rho \left( \bp^b, \Delta^b \right) \right)
& \leq &
m  \sum_{b = 0}^B \left(
\rho \left( \bp, \wh\Delta^b \right) -
\rho \left( \bp^b, \wh\Delta^b \right)
+ 2 L\, d \right) \\
& = & m \sum_{b = 0}^B \left(
\rho \left( \bp, \wh\Delta^b \right)
- \rho \left( \bp^b, \wh\Delta^b \right)\right) + 2 m (B+1) L\, d~.
\end{eqnarray*}
The terms $2 m (B+1) L \, d = 2 n L \, d$ are the first two terms in the
upper bound of the theorem.
The remaining part is bounded by using the same slope inequality argument
as in the previous section (recall that $\wt{r}$ denotes a sub-gradient bounded
between $-K$ and $K$):
\[
m \sum_{b = 0}^B \left(
 \rho \left( \bp, \wh\Delta^b \right)
- \rho \left( \bp^b, \wh\Delta^b \right) \right)
\leq m \sum_{b = 0}^B \wt{r} \left( \bp^b, \wh\Delta^b \right) \, \cdot \, \left( \bp - \bp^b \right)~.
\]
Finally, we deal with the mixing with the uniform distribution:
\begin{eqnarray*}
m \sum_{b = 0}^B \wt{r} \left( \bp^b, \wh\Delta^b \right) \, \cdot \, \left( \bp - \bp^b \right)
& \leq &
(1-\gamma) m \sum_{b = 0}^B \wt{r} \left( \bp^b, \wh\Delta^b \right) \, \cdot \, \left( \bp - \wt{\bp}^b \right)
+ 2K\, \gamma m (B+1) \\*
& & \mbox{(since, by definition, $\bp^b = (1-\gamma) \wt{\bp}^b + \gamma \bu$)} \\
& \leq & (1-\gamma) m \left( \frac{\ln N}{\eta} + \frac{(B+1)K^2 \eta}{2} \right) + 2K\, \gamma m (B+1) \\*
& & \mbox{
(by (\ref{expertsthm}))}  \\
& \le & \frac{m \ln N}{\eta} + \frac{n K^2\eta}{2} + 2K\, \gamma n ~.
\end{eqnarray*}
The proof is concluded by the Hoeffding-Azuma inequality which,
as in (\ref{HAzconcludes}), gives the final term in the bound.
The union bound indicates that the obtained bound holds with
probability at least $1 - (B+2) \delta \geq 1 - (n/m + 1) \delta$.
\end{proof}

\section{Deterministic feedback signals depending on action--outcome pair}
\label{ab}

In this last section we explain how in the case of deterministic
feedback signals the forecaster of the previous section can be modified so
that the order of magnitude of the per-round regret improves to
$n^{-1/3}$. This relies on the linearity of $\rho$ in its second
argument. In the case of random feedback signals, $\rho$ may not be linear
and it is because of this fact that we needed to group rounds of size $m$.
If the feedback signals are
deterministic, such grouping is not needed and the rate $n^{-1/3}$ is
obtained as a trade-off between an exploration term ($\gamma$) and the
cost payed for estimating the feedback signals ($\sqrt{1/(\gamma n)}$).  This
rate of convergence has been shown to be optimal in \cite{CeLuSt06}
even in the Hannan-consistent case.  The key property is summarized in
the next technical lemma, whose proof is postponed to the appendix.

\begin{lemma}
\label{LmLin}
For every fixed $\bp$, the function $\rho(\bp,\cdot)$ is linear on $\cF$.
\end{lemma}

\begin{remark}
The fact that the forecaster does not need to group rounds in the case
of deterministic feedback signals has an interesting consequence.  It is easy
to see from the proofs of Proposition~\ref{prop:aa} and
Theorem~\ref{Th:DeterministicActionOutcome}, through the linearity
property stated above, that the results presented there are still
valid when the payoff function $r$ may change with time (even, when
the environment can set it). The definition of the regret is then
generalized as
\[
R_n = \max_{\bp} \min_{z_1^n : H(\cdot,\ol\bz_n) = H(\cdot,\ol\bq_n)}
\frac{1}{n} \sum_{t=1}^n r_t(\bp,z_t) - \frac{1}{n} \sum_{t=1}^n r_t(I_t,J_t)~,
\]
where $\ol\bz_n$ is the empirical distribution of the sequence of
outcomes $z_1^n = (z_1,\ldots,z_n)$, and the same bounds hold. This
may model some more complex situations, including
Markov decision processes.  Note that choosing time-varying reward
functions was not possible with the forecasters of
\cite{PiSc01,CeLuSt06}, since these relied on a crucial structural
assumption on the relation between $r$ and $h$.
\end{remark}

Next we describe the modified forecaster. Denote by $\cH$ the vector space
generated by $\cF \subset \R^{| \cS | N}$ and $\Pi$ the linear
operator which projects any element of $\R^{| \cS | N}$ onto $\cH$.
Since the $\rho(\bp,\cdot)$ are linear on $\cF$, we may extend them
linearly to $\cH$ (and with a slight abuse of notation we write
$\rho$ for the extension). As a consequence, the functions
$\rho(\bp, \Pi(\cdot))$ defined on $\R^{| \cS |N}$ are linear and coincide with the original definition on $\cF$.  We denote by
$\wt{r}$ a sub-gradient (i.e., for
all $\Delta \in \R^{\abs[\cS] N}$, $\wt{r}(\cdot,\Delta)$ is a
sub-gradient of $\rho(\cdot,\Pi(\Delta))$).

The sub-gradients are evaluated at the following points.
(Recall that since the feedback signals are deterministic,
$s_t = h(I_t, J_t)$.) For $t = 1, 2, \ldots$, let
\begin{equation}
\label{defPD-group3}
\wh{h}_t = \left[ \wh{h}_{i,t} \right]_{i=1,\ldots,N}
= \left[ \frac{\delta_{s_t}}{p_{i,t}} \IND_{I_t = i} \right]_{i=1,\ldots,N}~.
\end{equation}
The $\wh{h}_{i,t}$ estimate the feedback signals $H(i,J_t) = \delta_{h(i,J_t)}$
received by action $i$
at round $t$. They are still conditionally unbiased estimators of the $h(i,J_t)$, and so is $\wh{h}_t$ for $H(\cdot,J_t)$.
The proposed forecaster is defined in Figure \ref{figab} and
the regret bound is established in Theorem~\ref{Th:DeterministicActionOutcome}.
\begin{figure}[t]
\bookbox{
\textbf{Parameters:} Real numbers $\eta, \gamma > 0$.
\\
\textbf{Initialization:} $\bw_1 = (1,\ldots,1)$.

\smallskip\noindent
For each round $t=1,2,\ldots$

\smallskip\noindent
\begin{enumerate}
\item choose the distribution $\bp_t = (1-\gamma) \wt\bp_t + \gamma \bu$,
where $\wt\bp_t$ is defined component-wise as
\[
\wt{p}_{k,t} = \frac{w_{k,t}}{\sum_{j=1}^N w_{j,t}}
\]
and $\bu$ denotes the uniform distribution, $\bu = (1/N,\ldots,1/N)$;
then draw an action $I_t$ from $\{1,\ldots,N\}$ according to $\bp_t$;
\item
perform the update
\[
    w_{k,t+1} = w_{k,t}\,e^{\eta\,\left(\wt{r}\left(\bp_t,\wh{h}_t \right)\right)_k}
\qquad
    \mbox{\rm for each $k=1,\ldots,N$,}
\]
where $\Pi$ is the projection operator defined after the statement of Lemma~\ref{LmLin},
for all $\Delta \in \R^{\abs[\cS] N}$,
$\wt{r}(\cdot,\Delta)$ is a sub-gradient of $\rho(\cdot,\Pi(\Delta))$,
and $\wh{h}_t$ is defined in (\ref{defPD-group3}).
\end{enumerate}
}
\caption{\label{figab}
The forecaster for deterministic feedback signals depending on action--outcome pair.
}
\end{figure}

\begin{theorem}
\label{Th:DeterministicActionOutcome}
There exists a constant $C$ only depending on $r$ and $h$ such that
for all $\delta > 0$, $\gamma \in (0,1)$, and
$\eta > 0$,
the regret for any strategy of the environment
is bounded, with probability at least $1 - \delta$, as
\begin{eqnarray*}
R_n & \le &
2 N C \sqrt{\frac{2}{n\gamma} \ln \frac{2}{\delta}} + \frac{2}{3} \frac{N C}{\gamma
n} \ln \frac{2}{\delta}
+ \frac{\ln N}{\eta n} + \frac{\eta K^2}{2}
+ 2K\,\gamma +
\sqrt{\frac{1}{2n} \ln \frac{2}{\delta}}~.
\end{eqnarray*}
The choice $\gamma \sim n^{-1/3} N^{2/3}$ and $\eta \sim \sqrt{(\ln N)/n}$
ensures that, with probability at least $1 - \delta$,
$R_n = O\left(n^{-1/3}  N^{2/3} \sqrt{\ln (1/\delta)} \right)$.
\end{theorem}

Note that here, as in Section~\ref{aa} (see Remark~\ref{r:eta}),
the tuning of the parameters can be done efficiently
on-line without resorting to the
``doubling trick.''
The optimization of
the upper bound (in both $\gamma$ and $\eta$) requires the knowledge
of $N$, $C$, $K$, and $n$. The first three parameters only depend on the
game and are known or may be calculated beforehand (the proof
indicates an explicit expression for $C$
and the bound on the sub-gradients may be computed as explained
in Section~\ref{subgrad}).
If $n$ and/or $K$ are unknown, their tuning
may be dealt with by taking time-dependent $\gamma_t$ and $\eta_t$.

\begin{proof}
The proof is similar to the one of Theorem~\ref{Th:RandomActionOutcome},
except that we do not have to consider the grouping steps
and that we do not apply the Hoeffding-Azuma inequality
to the estimated feedback signals
but to the estimated rewards. By the  bound on $r$ in terms of $\rho$
and convexity (linearity) of $\rho$ in its second argument,
\[
n \, \rho(\bp,H(\cdot,\ol\bq_n)) - \sum_{t=1}^n r(\bp_t,J_t)
 \leq
\sum_{t=1}^n \left( \rho \left( \bp, H(\cdot,J_t) \right)
                  - \rho \left( \bp_t, H(\cdot,J_t) \right) \right)~.
\]
Next we estimate
\[
\rho \left( \bp, H(\cdot,J_t) \right) - \rho \left( \bp_t, H(\cdot,J_t) \right) \qquad \mbox{by} \qquad
\rho \left( \bp, \Pi \left( \wh{h}_t \right) \right) - \rho \left( \bp_t, \Pi \left( \wh{h}_t \right) \right)~.
\]
By Freedman's inequality (see, again, \cite[Lemma A.1]{CeLuSt06}), since
$\wh{h}_t$ is a conditionally unbiased estimator of $H(\cdot,J_t)$ and all functions
at hand are linear in their second argument, we get that,
with probability at least $1 - \delta/2$,
\begin{eqnarray*}
\lefteqn{\sum_{t=1}^n \left( \rho \left( \bp, H(\cdot,J_t) \right)
- \rho \left( \bp_t, H(\cdot,J_t) \right) \right)} \\
& = & \sum_{t=1}^n \left( \rho \left( \bp, \Pi \left( H(\cdot,J_t) \right) \right)
- \rho \left( \bp_t, \Pi \left( H(\cdot,J_t) \right) \right) \right) \\*
& \leq &
\sum_{t=1}^n \left( \rho \left( \bp, \Pi \left( \wh{h}_t \right) \right)
- \rho \left( \bp_t, \Pi \left( \wh{h}_t \right) \right) \right)
+  2 N C \sqrt{2 \frac{n}{\gamma} \ln \frac{2}{\delta}}
+ \frac{2}{3} \frac{N C}{\gamma} \ln \frac{2}{\delta}
\end{eqnarray*}
where,
denoting by $\be_i(\delta_{h(i,j)})$ the column vector whose $i$-th component
is $\delta_{h(i,j)}$ and all other components equal 0,
\[
C = \max_{i,j} \max_{\bp} \rho \left( \bp, \Pi \left[ \be_i(\delta_{h(i,j)})
 \right] \right) < +\infty~.
\]
(A more precise look at the definition of $C$ shows that it is less than the
maximal $\ell_1$ norm
of the barycentric coordinates of the points
$\Pi [\be_i(\delta_{h(i,j)})]$
with respect to the $h(\cdot,j)$.)
This is because for all $t$, the conditional variances are bounded as follows. For all $\bp'$,
\begin{eqnarray*}
\EXP_t \left[ \rho \left( \bp', \Pi \left( \wh{h}_t \right) \right)^2 \right]
& = & \sum_{i=1}^N p_{i,t} \, \rho \left( \bp', \Pi \left[
\be_i(\delta_{h(i,j)}/p_{i,t})
\right] \right)^2 \\
& = & \sum_{i=1}^N \frac{1}{p_{i,t}} \, \rho \left( \bp', \Pi \left[
\be_i(\delta_{h(i,j)}/p_{i,t})
\right] \right)^2
 \leq  \sum_{i=1}^N \frac{C^2}{p_{i,t}} \leq \frac{C^2 N^2}{\gamma}~.
\end{eqnarray*}

The remaining part is bounded by using the same slope inequality argument
as in the previous sections
(recall that $\wt{r}$ denotes a sub-gradient in the first argument of $\rho(\cdot,\Pi(\cdot))$,
bounded between $-K$ and $K$),
\[
\sum_{t=1}^n \left( \rho \left( \bp, \Pi \left( \wh{h}_t \right)  \right)
- \rho \left( \bp_t, \Pi \left( \wh{h}_t \right)  \right) \right)
\leq \sum_{t=1}^n \wt{r} \left( \bp_t, \wh{h}_t \right) \, \cdot \, \left( \bp - \bp_t \right)~.
\]
Finally, we deal with the mixing with the uniform distribution,
\begin{eqnarray*}
\sum_{t=1}^n \wt{r} \left( \bp, \wh{h}_t \right) \, \cdot \, \left( \bp - \bp \right)
& \leq & (1-\gamma) \sum_{t=1}^n \wt{r} \left( \bp_t, \wh{h}_t \right) \, \cdot \, \left( \bp - \wt{\bp}_t  \right)
+ 2K\,\gamma n \\*
& & \mbox{(since by definition $\bp_t = (1-\gamma) \wt{\bp}_t + \gamma \bu$)} \\
& \leq & (1-\gamma)\left( \frac{\ln N}{\eta} + \frac{n\eta K^2}{2} \right)
 + 2K\,\gamma n
\quad \mbox{
(by (\ref{expertsthm})).}
\end{eqnarray*}
As before, the proof is concluded by the Hoeffding-Azuma inequality
(\ref{HAzconcludes}) and the union bound.
\end{proof}

\acknowledgments{We thank Tristan Tomala for helpful discussions.
Shie Mannor was partially supported by the Canada Research
Chairs Program and by the Natural Sciences and Engineering Research Council of Canada.
G\'abor Lugosi acknowledges the support of the Spanish Ministry of Science and Technology
grant MTM2006-05650 and
of Fundaci\'on BBVA. Gilles Stoltz was partially supported by the French ``Agence Nationale pour la Recherche''
under grant JCJC06-137444 ``From applications to theory in learning and adaptive statistics.''
G\'abor Lugosi and Gilles Stoltz acknowledge the PASCAL Network of Excellence under EC grant no.\ 506778.

An extended abstract of this paper appeared in the
\emph{Proceedings of the 20th Annual Conference on Learning Theory}, Springer, 2007.
}

\appendix

\section{Uniform Lipschitzness of $\rho$}

\begin{proposition} The function $(\bp,\Delta) \mapsto \rho(\bp,\Delta)$ is uniformly Lipschitz in its second argument.
\label{prop:uniformLip}
\end{proposition}
\begin{proof}
We consider the general case where the signal distribution depends on both the actions
and outcomes. Accordingly, we can write $\rho(\bp,\Delta)$ as the solution of the following linear
program (we denote $\Delta = (\Delta_1,\ldots,\Delta_N) \in \cF \subset \cP(S)^N$,
where, as usual, we identify each $\Delta_j$ with
a $|\cS|$-dimensional vector):
\[
\begin{array}{ccrl}
\rho(\bp,\Delta) & = & \min\limits_{\bq \in \R^{[\cS|}} & r(\bp,\cdot)^\top \bq \\
\\
& \mbox{s.t.\ } & H^k \, \bq & = \,\, \Delta_k~, \quad k=1,2,\ldots,N~, \\
& & \be_M^\top \, \bq & = \,\, 1~,\\
& & \bq & \ge \,\, 0~,
\end{array}
\]
where $r(\bp,\cdot) = (r(\bp,j))_j$ is an $M$-dimensional vector, $\be_M$ is an $M$-dimensional
vector of ones, and $H^k = H(k,\cdot)$ is the $|\cS| \times M$ matrix, whose entry $(s,j)$ is the probability of observing
signal $s$ when action $k$ is chosen and the outcome is $j$.

The program is feasible for every $\Delta \in \cal F$ so by the duality
theorem,
\begin{equation}
\label{eq:optproblem}
\begin{array}{ccll}
\rho(\bp,\Delta) & = & \max\limits_{\by \in \R^{N |\cS| +1}} \ \ \left[\Delta_1^\top \, \Delta_2^\top \,\ldots \,\Delta_N^\top \,1\right] \, \by \\
\\
& \mbox{s.t.\ } &  \left[H^1(\cdot,j)^\top \, H^2(\cdot,j)^\top \, \ldots \, H^N(\cdot,j)^\top \, 1 \right] \, \by & \leq \,\, r(\bp,j)~, \quad j=1,2,\ldots,M~, \\
& & \hfill \by & \geq \,\, 0~,
\end{array}
\end{equation}
where we recall that $H^k(\cdot,j)$ is the $|\S|$-dimensional vector whose $s$-th entry is the probability
of observing signal $s$ if the action is $k$ and the outcome is $j$.

We first claim that $\Delta \mapsto \rho(\bp,\Delta)$ is Lipschitz for every fixed $\bp$.
Indeed, for every fixed $\bp$ the optimization problem involves $\Delta$ only
through the objective function. We thus have that the solution to the optimization
problem is obtained at one of finitely many values of $\by$ (the vertices of the feasible
cone defined by the constraints of program (\ref{eq:optproblem})).
(More precisely, the obtained cone may be unbounded if there are
some unconstrained components of $\by$.
This happens when there exists an $s$ such that $H^k(s,j) = 0$ for all $j$.
But then $\Delta_k(s) = 0$ as well and we do not care about the unbounded component
$(k-1)N+s$ of $\by$.)
Since $\rho(\bp,\cdot)$ is a maximum of finitely many linear functions we
obtain that it is Lipschitz, with Lipschitz constant bounded by the
maximal $\ell_1$ norm of the vertices of the feasible
cone of (\ref{eq:optproblem}).

We now prove that the Lipschitz constant is uniform with respect to $\bp$.
It suffices to consider the polytope defined by
\[
\left\{ \by \in \R^{N |\cS| +1}: \,\, \by \geq 0, \ \left[
H^1(\cdot,j)^\top \, H^2(\cdot,j)^\top \, \ldots \, H^N(\cdot,j)^\top \, 1 \right] \, \by \leq 1, \quad j=1,2,\ldots,M
 \right\}~.
\]
This is a cone, and the vertex $\by$ with the maximum $\ell_1$ norm upper bounds
the Lipschitz constant of the $\rho(\bp,\cdot)$, for all $\bp$. (As before, any
unbounded components of $\by$ do not matter to the optimization problem.)
\end{proof}

\begin{remark}
\label{rem:boundL}
{\rm Observe from the proof that an upper bound on the uniform Lipschitz constant can
be easily computed by solving the following linear program,
\begin{equation*}
\begin{array}{ccll}
  & &\max \limits_{\by \in \R^{N |\cS| +1}} \ \
 \be_{NS+1}^\top\, \by \\
\\
& \mbox{s.t.\ } &  \left[H^1(\cdot,j)^\top \, H^2(\cdot,j)^\top \,
\ldots \, H^N(\cdot,j)^\top \, 1 \right] \, \by & \leq \,\, 1~, \quad j=1,2,\ldots,M~, \\
& & \hfill \by & \geq \,\, 0~.
\end{array}
\end{equation*}
}
\end{remark}

\section{Proof of Lemma~\ref{LmLin}}

It is equivalent to prove that for all fixed $\bp$, the function $\bq \mapsto \rho(\bp,H(\cdot,\bq))$ is linear
on the simplex.
Actually, the proof exhibits a simpler expression for $\rho$.

To this end, we first group together the outcomes with same feedback signals and
define a mapping
\[
T : \cP \bigl( \{1,\ldots,M\} \bigr) \to \cP \bigl( \{1,\ldots,M\} \bigr)~,
\]
where $\cP(\{1,\ldots,M\})$ is the set of all probability distributions $\bq$ on
the outcomes. Formally, consider the
binary relation
defined by $j \equiv j'$ if and only if $h(\cdot,j) =
h(\cdot,j')$. (We use here the notation $h$ to emphasize that we deal
with deterministic feedback signals.) Denote by $F_1,\ldots,F_{M'}$ the
partition of the outcomes $\{1,\ldots,M\}$ obtained so, and pick
in every $F_j$ the outcome $y_j$ with minimal reward $r(\bp,y_j)$ against $\bp$
(ties can be broken arbitrarily, e.g., by choosing the outcome
with lowest index).
Then, for every $\bq$, the distribution $\bq' = T(\bq)$ is defined as
$q'_{y_j} = \sum_{y \in F_j} q_y$, for $j= 1,\ldots,M'$, and
$q'_k = 0$ if $k \ne y_j$ for all $j$.

$T$ is a linear projection (i.e., $T \circ T = T$).
It is easy to see that in the case of deterministic feedback signals,
$H(\cdot,\bq) = H(\cdot,\bq')$ if and only if $T(\bq) = T(\bq')$.
This implies that
\begin{equation}
\label{exprerho}
\rho(\bp,H(\cdot,\bq)) = \min_{\bq' \, : \, T(\bq') = T(\bq)} r(\bp,\bq') = r(\bp,T(\bq))
\end{equation}
where the last equality follows from the fact that, by choices of the $y_j$,
$r(\bp,\bq') \geq r(\bp,T(\bq'))$ for all $\bq'$, with equality
for $\bq' = T(\bq) = T^2(\bq)$.
By linearity of $T$,
$\ \bq \mapsto r(\bp,T(\bq)) = \rho(\bp, H(\cdot,\bq))$ is therefore linear
itself, as claimed.

Note that the equivalence of $H(\cdot,\bq) = H(\cdot,\bq')$ and
 $T(\bq) = T(\bq')$,
together with (\ref{exprerho}),
implies the following
sufficient condition for Hannan-consistency (for necessary and sufficient conditions, see
\cite{PiSc01,CeLuSt06}). It is more general than the distinguishing actions condition
of \cite{CeLuSt06}.
\begin{remark}
\label{rkHC}
Whenever $H$ has no two identical columns in the case
of deterministic feedback, i.e., $h(\cdot,j) \ne h(\cdot,j')$ for all $j \ne j'$,
one has that for all $\bp$ and $\bq$,
\[
\rho(\bp,H(\cdot,\bq)) = r(\bp,\bq)~.
\]
\end{remark}
The condition is satisfied, for instance, for multi-armed bandit problems, where $h = r$ (provided that we identify outcomes
yielding the same rewards against all decision-maker's actions).


\end{document}